\documentclass{amsart}

\usepackage{enumerate}
\usepackage{ifpdf}
\usepackage{amsmath}
\usepackage{amsfonts}   % if you want the fonts
\usepackage{amssymb}    % if you want extra symbols
\usepackage{amsthm}
\usepackage{hyperref}
\usepackage{setspace}

\newcommand{\ignore}[1]{}

\renewcommand{\Re}{\operatorname{Re}}
\renewcommand{\Im}{\operatorname{Im}}
\newcommand{\Orb}{\operatorname{Orb}}
\newcommand{\hol}{\operatorname{hol}}

\newcommand{\abs}[1]{\left\lvert {#1} \right\rvert}

\newcommand{\C}{{\mathbb{C}}}
\newcommand{\R}{{\mathbb{R}}}

\newcommand{\sC}{{\mathcal{C}}}

\newcommand{\sH}{{\mathcal{H}}}
\newcommand{\sM}{{\mathcal{M}}}

\newcommand{\sU}{{\mathcal{U}}}

\newcommand{\sX}{{\mathcal{X}}}

\newcommand{\sZ}{{\mathcal{Z}}}
\newcommand{\fS}{{\mathfrak{S}}}

\newcommand{\dhat}[1]{\hat{\hat{#1}}}

\newtheorem{thm}{Theorem}[section]
\newtheorem*{thmnonum}{Theorem}

\newtheorem{prop}[thm]{Proposition}
\newtheorem{cor}[thm]{Corollary}
\newtheorem{lemma}[thm]{Lemma}

\theoremstyle{definition}
\newtheorem{defn}[thm]{Definition}

\theoremstyle{remark}

\author{Ji\v{r}\'i Lebl}
\address{ Department of Mathematics, University of California
at San Diego, La Jolla, CA 92093-0112, USA}
\email{jlebl@math.ucsd.edu}
\date{February 24, 2007}

\ifpdf
\fi

\title[Nowhere minimal CR submanifolds and Levi-flat
hypersurfaces]{Nowhere minimal CR submanifolds and Levi-flat
hypersurfaces}

\begin{document}

%\doublespace

\begin{abstract}
A local uniqueness property of holomorphic functions on real-analytic
nowhere minimal CR submanifolds of higher
codimension is investigated.  A sufficient condition called
almost minimality is given and studied.
A weaker necessary condition, being contained
a possibly singular real-analytic Levi-flat hypersurface is studied and
characterized.  This question is completely
resolved for algebraic submanifolds of codimension 2 and
a sufficient condition for noncontainment is given for non algebraic
submanifolds.  As a consequence, an example of a submanifold of codimension
2, not biholomorphically equivalent to an algebraic one, is given.
We also investigate the structure of singularities of Levi-flat
hypersurfaces.
\end{abstract}

\maketitle

%\enlargethispage{\baselineskip}

\section{Introduction}

In this paper we investigate some local properties of
nowhere minimal real-analytic
CR submanifolds of higher codimension.
In particular we are interested in a modulus uniqueness property for holomorphic
functions, that is, when is a holomorphic function uniquely determined (up to
a unimodular constant) by
its modulus on a CR submanifold.  We introduce a sufficient geometric
condition called almost minimality, and we study related properties
of such submanifolds.  We also introduce
a necessary condition, that is, being contained in a singular Levi-flat
hypersurface, and thus we will find it necessary to study the structure
of the singular set of such hypersurfaces.

Background material is taken mostly from \cite{BER:book}.  We first fix some
terminology.
Let $M \subset \C^N$ be real-analytic submanifold defined near the origin.
The tangent vectors of the form
$\sum_{j=1}^N a_j \frac{\partial}{\partial \bar{z}_j}$ tangent to $M$
are called the
\emph{CR vectors}.  If this space has constant dimension on $M$,
the submanifold is said to be a \emph{CR submanifold}, and
the complex dimension of the CR tangent space is called
the \emph{CR dimension} of $M$.  If a
CR submanifold is not contained in a proper complex analytic subvariety,
we say it is a \emph{generic submanifold}.
We denote by $\Orb_p$ the \emph{CR orbit} of $M$ at $p$,
that is, the germ of a CR submanifold of $M$ through $p$ of smallest
dimension that has the same CR dimension as $M$.
If $\Orb_p = M$ as germs, then
$M$ is said to be minimal at $p$.  If a real-analytic submanifold is
minimal at one point, then it is minimal outside a real-analytic subvariety.
If $M$ is contained in a real-algebraic subvariety of $\C^N$ of the same
dimension as $M$, then $M$ is said to be \emph{real-algebraic}.
We will say that a generic submanifold is \emph{Levi-flat}, if there exist local
holomorphic coordinates $z=(z_1,\dots,z_N)$, such that
$M$ can be given by $\Im z_1 = \cdots = \Im z_d = 0$.
A real-analytic, possibly singular hypersurface $H$ (defined by the vanishing
of a single real-analytic real valued function) is said to be
\emph{Levi-flat},
if near the nonsingular points of hypersurface dimension,
there exist local holomorphic coordinates,
such that $\Im z_1 = 0$ defines $H$.  That is, Levi-flat hypersurfaces are
locally (near nonsingular points) foliated by complex analytic hypersurfaces,
and we call this the
\emph{Levi foliation}.
We denote by $H^*$
the nonsingular points of $H$ that are of hypersurface dimension, that is,
the points near which $H$ is a real-analytic submanifold of real codimension
1.  We
define $H_s := H \backslash H^*$.

All minimal submanifolds have the modulus uniqueness property.
The story is not so simple with nowhere minimal
submanifolds.  We introduce a sufficient geometric condition,
called almost minimality, for
a submanifold to have the modulus uniqueness property.
For a generic submanifold $M$ through the origin,
$M$ is almost minimal at 0, if for any connected neighbourhood $U$ of 0,
there exists a point $p \in U$,
such that the CR orbit at $p$
is not contained in a proper complex analytic subvariety of $U$.
For example (see \S \ref{sec:ex}) the manifold $M_\lambda$,
given by the defining equations
\begin{align*}
\bar{w}_1 &= e^{iz\bar{z}} w_1 ,
\\
\bar{w}_2 &= e^{i\lambda z\bar{z}} w_2 ,
\end{align*}
is almost minimal at 0 if and only if $\lambda$ is irrational.

The modulus
uniqueness property described above is equivalent to the submanifold being
locally contained in a
possibly singular real-analytic Levi-flat hypersurface defined by
the vanishing of the imaginary part of a meromorphic function.  See
\S \ref{sec:uniq} for discussion of the modulus uniqueness property.
We will
consider a weaker condition, that is, when the
submanifold is contained in any Levi-flat hypersurface $H$.
%If we are in $\C^N$,
%we consider a real-analytic subvariety $H$ of real dimension
%$2N-1$ Levi-flat if it is Levi-flat at all the nonsingular points of highest
%dimension, which we denote by $H^*$ and we define
%$H_s := H \backslash H^*$.
In particular,
we will be concerned about
when our higher codimension manifold $M$ is contained in
$\overline{H^*}$, which is not necessarily the same as $H$.  Singular Levi-flat
hypersurfaces with quadratic tangent cones have been studied by Burns and
Gong in \cite{burnsgong:flat}, and a similar approach,
studying the Segre varieties of $H$, is taken in this paper.  Burns and
Gong also give an example of a Levi-flat hypersurface not defined by the
vanishing of the real part of a meromorphic function.  So being contained in
a Levi-flat hypersurface is a potentially weaker condition on $M$ than
having the modulus uniqueness property.  Singular
Levi-flat hypersurfaces have also been studied by Bedford
\cite{bedford:flat} in case
the singularity is contained in a codimension 2 complex variety.

For a generic submanifold $M \subset \C^N$, we will consider $M$ in normal
coordinates $(z,w) \in \C^n \times \C^d$, where
$d$ is the real codimension of $M$ and $n$ is the CR dimension of $M$,
and
$M$ is given by $w = Q(z,\bar{z},\bar{w})$, where $Q$ is a holomorphic
mapping defined
in a neighbourhood of the origin in $\C^n \times \C^n \times \C^d$,
$Q(0,\zeta,\omega) \equiv Q(z,0,\omega) \equiv \omega$, and
$Q(z,\zeta,\bar{Q}(\zeta,z,w)) \equiv w$.  We should note, however, that
normal coordinates are not unique.  The first main result of this paper is
the following.

\begin{thm} \label{fullprojwflatthm}
% Original labels of the following are:
%label projwflatthm
%label algcodim2thm
Let $M$ be a germ of a
generic real-analytic
codimension 2 submanifold through the origin
given in normal coordinates $(z,w)$ and let $M$ be nowhere minimal.  
Then $M \subset \overline{H^*}$, where $H$ is a germ of
a possibly singular real-analytic
Levi-flat hypersurface if and only if 
the projection of $M$ onto the second factor in $(z,w)$
is contained in a germ of a
possibly singular real-analytic hypersurface.
If $M$ is real-algebraic, then such an $H$ always exists and is real-algebraic.
Moreover, if $M$ is not Levi-flat, then $H$ is unique.
\end{thm}

Thus to answer the question of when $M$ sits inside a Levi-flat hypersurface,
it is sufficient to study the projection of
$M$ onto the second factor in normal coordinates, if $M$ is not
real-algebraic.
The question is fully answered in case
$M$ is real-algebraic.
In fact, if such an $H$ exists,
then there exists one defined by an equation independent of $z$.
Further, we will show that if two holomorphic functions $f$ and $g$
have equal modulus on $M$, then $f/g$ also only depends on $w$.

We also prove that an almost minimal submanifold of codimension 2 cannot be
contained in a Levi-flat hypersurface.  So an almost minimal submanifold
that is nowhere minimal is an example of a
submanifold, which is not locally biholomorphic to a real-algebraic
submanifold.  For $\lambda$ irrational, the submanifolds $M_\lambda$
defined above are such examples.
Examples of such
hypersurfaces both minimal and nonminimal can be found in
\cite{BER:localprop} and \cite{HJY:example}.

When $M$ is a almost minimal at $p$,
we will study 
the dimension of $\hol(M,p)$, the space of
\emph{infinitesimal holomorphisms}
at $p$, that is, the Lie algebra generated by germs at $p$
of real-analytic vector fields $X$ on $M$ defined in some neighbourhood $U$
of $p$, such that for each $q \in U$ there is another neighbourhood $q \in V
\subset U$ such that the map $z \mapsto \exp tX \cdot z$ for $\abs{t} \leq
\epsilon$ is a CR diffeomorphism of $M$ (a diffeomorphism of $M$ that
preserves the CR vector bundle).

A vector field $X$ in $\C^N$
is called a \emph{holomorphic vector field}, if we can write it locally as
$X = \sum_{k=1}^N a_k(z) \frac{\partial}{\partial z_k}$, where the $a_k$
are holomorphic in $z \in \C^N$.
A submanifold $M$
is said to be \emph{holomorphically nondegenerate} at $p \in M$,
if there does not exist
any germ at $p$ of a nonzero holomorphic vector field tangent to $M$.
If $M$ is connected, real-analytic and generic it turns out, that if it is
holomorphically nondegenerate at one point it is so at all points.
Being holomorphically nondegenerate is
a necessary condition for $\dim_\R \hol(M,p) < \infty$.
In the case $M$ is a hypersurface
Staton \cite{Stanton:crauthyp} proved that this is in fact a sufficient
condition.  For higher codimension submanifolds,
Baouendi, Ebenfelt and Rothschild \cite{BER:craut} proved that 
$\dim_\R \hol(M,p) < \infty$ if $M$ is minimal at $p$, and if $M$
is nowhere minimal, then on a dense open subset of $M$,
$\dim_\R \hol(M,p)$ is either zero or infinite.  We prove the following
result for almost minimal submanifolds.

\begin{thm}
\label{holfindim}
Let $M \subset \C^N$ be a connected, real-analytic holomorphically
nondegenerate generic submanifold and suppose
$p \in M$ and $M$ is almost minimal at $p$.  Then
\begin{equation*}
\dim_\R \hol(M,p) < \infty .
\end{equation*}
\end{thm}

Finally, it will be necessary to know something about the structure of
the singular set of a Levi-flat hypersurface.  This result is also
of interest on its own.  We prove a technical theorem
in \S \ref{sec:lfsing} which has the following corollary.

\begin{thm} \label{flatopendense}
Let $H \subset \C^N$ be a singular
real-analytic
Levi-flat
hypersurface, and
let $M \subset H_s \cap \overline{H^*}$
be a smooth submanifold.
Then for $p$ on an open dense set of $M$, the germ of $M$ at $p$ is
contained in some germ of a complex variety or generic real-analytic
Levi-flat submanifold 
of real dimension $2N-2$.
\end{thm}

The paper has the following organization.  In \S \ref{sec:lfsing} we
study the singular set of Levi-flat hypersurfaces and prove
Theorem \ref{flatopendense}.  In \S \ref{sec:uniq} we study the
modulus uniqueness property.  In \S \ref{sec:flatness} we consider
when $M$ is contained in a Levi-flat hypersurface and prove
the first part of Theorem \ref{fullprojwflatthm}.
In \S \ref{sec:almostmin} we define
and study the almost minimality condition and in \S \ref{sec:craut}
we prove Theorem \ref{holfindim}.
In \S \ref{sec:alg} we study real-algebraic submanifolds and prove
the remainder of
Theorem \ref{fullprojwflatthm}.  Finally in \S \ref{sec:ex} we work out the
example $M_\lambda$ family of submanifolds and prove a slightly
more general result which can be used for generating further
examples of almost minimal submanifolds.

The author would like to thank Peter Ebenfelt for fruitful discussions and
guidance in the preparation of these results and of this
manuscript.  Also, the author would like to acknowledge M. Salah Baouendi,
Nordine Mir and Linda Preiss Rothschild for reading over the manuscript
and their useful comments.  Finally, the author would like to thank the
referee for corrections and improvements to the exposition.

\section{Singularity of Levi-flat hypersurfaces} \label{sec:lfsing}

As we noted before we will consider a singular real-analytic
Levi-flat hypersurface through the origin $H \subset U$ where
$U$ is an open neighbourhood of the origin in $\C^N$ given by
a $H = \{ z \in U \mid \rho(z,\bar{z}) = 0 \}$ for a real valued real-analytic
function $\rho$.  As we are interested in local properties
of $H$ we will assume that $U$ is small enough such that
$\rho$ can be complexified to $U \times {}^* U$ where ${}^* U = \{ z \mid
\bar{z} \in U \}$.  Further, we will assume that $U$ is connected.
As before we will denote by $H^*$
the nonsingular points of dimension $2N-1$.
Then we let $H_s := H \backslash H^*$.
We note that it is not necessarily true that $\overline{H^*} = H$, even
if $H$ is irreducible.
We say that $H$ is Levi-flat, if near each $p \in H^*$ there are
suitable holomorphic coordinates such that $H$ is given by
$\Im z_1 = 0$.  Burns and Gong \cite{burnsgong:flat} prove the
following useful lemma.

\begin{lemma} \label{flatsomewhereflateverywhere}
Let $H$ be an irreducible singular real-analytic hypersurface.  Then
if $H$ is Levi-flat at a single point of $H^*$, then it is
Levi-flat at all points of $H^*$.
\end{lemma}

Our main result about Levi flat hypersurfaces is the following 
theorem.

\begin{thm} \label{flathypersing}
Let $H \subset U \subset \C^N$ be a singular
real-analytic
Levi-flat
hypersurface.  Then
\begin{equation*}
H_s \cap \overline{H^*} \subset \bigcup_{j=1}^\infty M_j
\end{equation*}
where $M_j \subset U_j$ for some countable collection of open sets
$U_j \subset U$, and where $M_j$ is either a proper complex analytic
subvariety of $U_j$
or a generic real-analytic Levi-flat submanifold of real dimension at most $2N-2$.
\end{thm}

%Where we say that a generic submanifold is Levi-flat if in suitable
%coordinates it can be given by $\Im z_1 = \dots = \Im z_d = 0$.
Theorem \ref{flatopendense} in the
introduction follows from this technical result.

\begin{proof}[Proof of Theorem \ref{flatopendense}]
By the above theorem,
$M \subset \bigcup M_j$.  Suppose that there is no point in $M$ such that
near that point $M \subset M_j$ (as germs) for some $j$.
That means, $M \cap M_j$
is nowhere dense in $M$ (it does not contain an open set).  But there are
only countably many such sets, and so by Baire category theorem they cannot
cover all of $M$, which would be a contradiction.  Thus there has to exist a
point $p$ where $M$ is contained (as a germ at $p$) in some $M_j$.
This holds on an open set near $p$ as well, and furthermore,
since it holds for all open $U$, by taking $U$ smaller we can see that it has
to hold on an open dense set of $M$.
\end{proof}

A useful weaker result, at a point where $H_s$ is a submanifold of codimension
one in $H$, is the following.

\begin{cor} \label{flathypersingcor}
Let $H$ be a singular real-analytic Levi-flat hypersurface
defined in a neighbourhood
of the origin in $\C^N$, and
suppose that $H_s$ is a manifold
of dimension $2N-2$ and $H_s \subset \overline{H^*}$.  Then
$H_s$ is either complex analytic or Levi-flat.
\end{cor}

\begin{proof}
If $H_s$ was of a different type, then all the Levi-flat and complex analytic
$M_j$'s
have an intersection of a lower dimension with $H_s$.  By Baire
category theorem again, this is not possible, as there are only countably many.
\end{proof}

Thus we have a complete categorization of singularities if they are of
highest possible dimension and are in the closure of the nonsingular points.
There are examples where the singular set is
complex
(e.g. $\{ z \mid \Im z_1^2 = 0 \}$)
or Levi-flat (e.g. $\{z \mid (\Im z_1)(\Im z_2) = 0 \}$),
but it is not clear that an irreducible hypersurface can have a Levi-flat
singularity.

A smooth CR submanifold is said to be of \emph{finite type} at $p \in M$
if the CR vector fields,
their complex conjugates and finitely many commutators, span the
complexified tangent space $\C T_pM$ ($\C T_pM = \C \otimes_\R T_pM$).
In case $M$ is real-analytic, being finite type at $p$
is equivalent to being minimal at $p$.  It is not hard to see that
if $M$ is finite type at $p$, then there cannot exist a holomorphic
function in a neighbourhood of $p$ which is real valued on $M$.
We can now rule out all smooth finite type
generic 
submanifolds of any codimension being contained in Levi-flat hypersurfaces.

\begin{cor} \label{minhnotminimal}
Let $H \subset U \subset \C^N$ be singular real-analytic
Levi-flat hypersurface, and let $M \subset \overline{H^*}$
be a
smooth
generic submanifold.
Then $M$ is not of finite type at any point.
\end{cor}

\begin{proof}
Take a point $p \in M$.
If $p \in M \cap H^*$, then there is some neighbourhood of $p$,
where in suitable local coordinates $H$ is given by $\Im z_1 = 0$, and thus
$z_1$ is real valued on an open set of $M$.
Since if $M$ would be of finite type at $p$, it would
be of finite type in a neighbourhood of $p$.  If there exists
a real valued holomorphic function on $M$ near $p$, $M$ cannot be
of finite type at $p$.
So let $p \in H_s$.  Again if $M$ would be of finite type at $p$ then 
it would be so near $p$, and there would either
be a point $q \in M \cap H^*$ where $M$ was of finite type, which we now
know cannot happen, or $M \subset H_s$ as germs at $p$.  But then by
Theorem
\ref{flatopendense} for some point $q \in M$, where $M$ would be
of finite type, it would be contained as germ
in either a complex variety or a Levi-flat generic submanifold which
is again impossible.  Thus $M$ cannot be of finite type.
\end{proof}

Before going into the proof of Theorem \ref{flathypersing},
let's fix some notation and background.
Let $\Sigma_z$ be the \emph{Segre variety} of $H$
at the
point $z$, that is the set
$\{ \zeta \in U \mid \rho(\zeta,\bar{z}) = 0 \}$,
and let $\Sigma_z'$ be the branches of $\Sigma_z$ completely inside $H$.
We say that $\Sigma_z$
is \emph{degenerate} if $\Sigma_z$ contains an open set of $\C^N$, that is, if
$\Sigma_z = U$ if $U$ is connected.

We will need some lemmas about Levi-flat hypersurfaces.  Both of the
following are given (in more generality) and proved in \cite{burnsgong:flat}.

\begin{lemma}
\label{irredlemma}
If $\rho$ is an irreducible germ of a real-analytic function near $0$
in $\C^N$, and $H := \{ z \mid \rho(z,\bar{z}) = 0 \}$ has dimension 
$2N-1$, then for any neighbourhood $U$ of 0, there is a smaller neighbourhood
$U' \subset U$ of 0, such that if $\hat{\rho}$ is any real-analytic function
on $U$
which vanishes on an open set of $H^* \cap U'$, then $\rho$ divides
$\hat{\rho}$ on $U'$.  Further, $\rho$ is irreducible as a germ of a holomorphic
function near origin in $\C^{2N}$.
\end{lemma}

\begin{lemma}
\label{sprimenonempty}
Let $H \subset U$ be as above and Levi-flat, and
suppose $z \in H$ is such that
$\Sigma_z$ is non-degenerate.  Then $\Sigma_z'$ is non-empty, and further one
branch of $\Sigma_z'$ passes through $z$.  If $z \in H^*$, then $\Sigma_z'$ has only
one branch through $z$, and this is the unique germ of a complex variety
through $z$.
\end{lemma}

Also, since we could pick $U$ smaller and smaller,
one branch of $\Sigma_z'$ must therefore always pass through $z$.

If $\rho$ is a defining function for
$H$ in a neighbourhood $U$, then
at all points of $H_s$, $\rho$ must have a vanishing gradient, since
otherwise $H$ would be a nonsingular hypersurface at that point.
In fact, picking a possibly smaller $U$,
$\{ z \in H \mid \partial \rho(z,\bar{z}) = 0 \}$ is a proper subvariety
of $H$ containing $H_s$
(here $\partial$ means the exterior derivative in the $z$
variables).  Assume $H$ is irreducible, complexify $\rho$ into 
$U \times {}^* U$, and let
$\sH = \{ (z,\zeta) \in U \times {}^* U \mid \rho (z,\zeta)
= 0 \}$.  Then by Lemma \ref{irredlemma}, $\rho$ is irreducible
as a holomorphic function (in a possibly smaller neighbourhood),
and thus generates the ideal of $\sH$ by the Nullstellensatz at every point
in $U \times {}^* U$.  Therefore,
the gradient of the complexified $\rho$ does not vanish at all nonsingular points
of $\sH$.  Near any $p \in H^*$ we have a local defining function with
nonvanishing gradient near
$p$, which when complexified divides $\rho$.  That means, 
near $p$, $H^*$ complexifies to a
germ of a smooth complex hypersurface in 
$U \times {}^* U$ contained in $\sH$.  Since $H^*$ is totally real
in this complex hypersurface we know $\partial \rho$ cannot vanish
identically on $H^*$ (or it would vanish in all
of $\sH$ since it is irreducible).  Hence, $\partial \rho = 0$ defines a
proper lower dimensional subvariety
of $H$ which contains $H_s$.  We can't quite say it equals $H_s$, as a point
$p$ could be in $H^*$, but the point $(p,\bar{p})$ could a priory be a singular
point of $\sH$.

\begin{lemma}
\label{interleviflat}
Let $H_1, H_2 \subset \C^N$ be
two connected nonsingular real-analytic Levi-flat
hypersurfaces, such that $0 \in H_1 \cap H_2$.
If $U$ is a sufficiently small neighbourhood
of $0$, and $H_1 \cap U \not= H_2 \cap U$,
then there exists a possibly empty proper complex analytic subvariety
$A \subset U$ such that
$(U \cap H_1 \cap H_2) \setminus A$ is either empty or a
generic real-analytic Levi-flat submanifold of codimension 2.
\end{lemma}

\begin{proof}
We let $U$ be small enough such that $H_1 \cap U$ and $H_2 \cap U$ are closed
in $U$ and hence we can assume that $H_1, H_2 \subset U$.  Further,
let $U$ be small enough such that there exist holomorphic coordinates in $U$
where
$H_1$ is given by $\Im z_1 = 0$ and $H_2$ is given by $\Im f = 0$,
where $f$ is holomorphic
with nonvanishing differential.  The set where the complex differentials of
$f$ and $z_1$ are linearly dependent is a complex analytic subvariety.
If the complex differentials are everywhere linearly dependent then
$f$ depends only on $z_1$ and thus the intersection of $H_1$ and $H_2$ is
complex analytic.  So suppose that outside a subvariety $A$, $f$
and $z_1$ have linearly independent differentials so locally in an even
smaller neighbourhood we can change coordinates again to make $f = z_2$ and
then the intersection is locally defined as $\Im z_1 = \Im z_2 = 0$ and
we are done.
\end{proof}

\begin{proof}[Proof of Theorem \ref{flathypersing}]
Recall that to prove the Theorem, we will cover $H_s \cap \overline{H^*}$ by 
countably many Levi-flat submanifolds of codimension 2 and local
complex analytic subvarieties.  These submanifolds and subvarieties need
not lie in $H$ itself, we just want their union as sets to contain
$H_s \cap \overline{H^*}$.

Let $H_s' := H_s \cap \overline{H^*}$.  The place in the proof where
we fail to cover all of $H_s$, if
$H_s \not\subset \overline{H^*}$, is in the application of
Lemma \ref{sprimenonempty}.

Assume that $H$ is irreducible.  If it is reducible, and we
prove the result for each branch, then it is also true for the union of
those branches.  This is because if $K$ and $L$ are
branches of $H = K \cup L$, then
$H_s = K_s \cup L_s \cup S$, where
$S$ is the set of points of $K^* \cap L^*$, where
$K^* \cap L^*$ is not a hypersurface.  Hence, if we have 
covered $K_s$ and $L_s$, the only other points that need to be covered
are points of $S$.  If $p \in S$ we pick a small enough neighbourhood
of $p$ and apply Lemma \ref{interleviflat}.
We can also assume $H$ it is irreducible in
arbitrarily small neighbourhoods of 0 as well for the same reason (so
irreducible as a germ).

First we note that the points $z \in U$ where $\Sigma_z$
is degenerate are inside a complex
analytic variety, because $z \in \Sigma_w$ implies by reality of $\rho$
that $w \in \Sigma_z$.  So that means that if $z$ is a degenerate point, then it
is contained in $\Sigma_w$ for all $w \in U$,
and thus is inside a complex analytic
subvariety $A$.  Because we only care about a countable union of local
varieties and manifolds, we can just cover
$U \setminus A$ by smaller neighbourhoods and work there.
Thus we can assume that $U$ contains no degenerate points.

Suppose $0 \in H$, and
suppose that a branch of $\Sigma_0$, call it $A$ again,
is contained in $H_s$.  Again, since we only care about a countable
union of local varieties and manifolds, we can cover $U \setminus A$
by small neighbourhoods and work there.  Thus we can assume that $\Sigma_0$
has no branch that is contained in $H_s$ (and thus not in $H_s'$).

By Lemma \ref{sprimenonempty}, $\Sigma_0'$ is non-empty and we now know that
no branch of it is contained completely in $H_s'$.  So we know
that there exists a point $\zeta \in \Sigma_0'$ such that $\zeta \in H^*$.
As $\Sigma_0'$ at $\zeta \in H^*$ is the unique complex variety
(again by Lemma \ref{sprimenonempty})
passing through $\zeta$ we know that $\Sigma_{\zeta}'$ shares this branch
with $\Sigma_0'$.

We can of course pick this $\zeta$ in a topological component of $(\Sigma_0')^*
\cap H^*$,
where $(\Sigma_0')^*$ is the nonsingular part of $\Sigma_0'$,
such that 0 is in the closure of this component.  As no branch of $\Sigma_0'$ lies
inside $H_s'$ and there is at least one branch through 0, then
at least one topological
component of $(\Sigma_0')^* \cap H^*$ will be such that $0$ is in its closure.

We look at a small neighbourhood $V$ of $\zeta$
such that $H \cap V$ is connected and nonsingular, and
further, such that $H$ is defined in $V$ by $\Im f(z) = 0$,
for some $f$ holomorphic in $V$ where the gradient of $f$ does not vanish
in $V$.

Pick a nonsingular
real-analytic curve $\gamma \colon (-\epsilon,\epsilon) \to H$
such that $\gamma(0) = \zeta$, $\{\gamma\} \subset V$, and furthermore,
that $\gamma$ is transverse to the Levi foliation of $H^*$.  We can do this
by just changing coordinates in $V$ such that $z_n = f$, and then our curve
might be $t \mapsto t \alpha$ where $\alpha \in C^n$ and $\alpha_n$ is not real.
Once we have $\gamma$ we can look
at the sets $\Sigma_{\gamma(t)}$ for various $t$.  These are
given by $\{ z \mid \rho(z,\overline{\gamma(t)})=0\}$.
However, we can just look at the zero set of the function
$(z,t) \mapsto \rho(z,\bar{\gamma}(t))$ as $t$ is real.
Further, we can pick $\gamma$ such that
$\rho(0,\bar{\gamma}(t))$ is not identically zero 
since if it were for all choices of $\gamma$ (by varying $\alpha$
above), then
$\Sigma_0$ would contain an open set in $H^*$ and thus 
would be degenerate, and we assumed it was not.
We can complexify $t$ and look at the zero set of
$\rho(z,\bar{\gamma}(t))$ in $U \times D_\epsilon$ (where $D_\epsilon$ is
the disk of radius $\epsilon > 0$).

Next apply the Weierstrass preparation theorem, which we can do in some
neighbourhood of $(0,0)$ in $U' \times D_{\epsilon'} \subset U \times
D_\epsilon$ and we get a polynomial
\begin{equation*}
F(z,t) = t^m + \sum_{j=0}^{m-1}a_j(z)t^j ,
\end{equation*}
whose zero set
is the zero set of $\rho(z,\bar{\gamma}(t))$.  Outside of the
discriminant set of $F$, $\Delta \subset U'$, we have (locally)
$m$ holomorphic functions $\{e_j\}_1^m$ which give us the
solutions to $F(z,e_j(z))=0$.  We look at
the places where these solutions are real, that is
the points in $U'$ where 
$e_j-\bar{e_j} = 0$.  To be able to complexify we look at the
function
\begin{equation*}
i^m \prod_{j,k=1}^m (e_j(z)-\overline{e_k(z)}).
\end{equation*}
It is easy to see
that this is a real function.  Furthermore, it is symmetric
both in the $e_j(z)$ and the $\overline{e_k(z)}$, this means that
after complexification we have
a well defined holomorphic function in
$(U' \times {}^*{U'}) \setminus (\Delta \times {}^*{\Delta})$,
and continuous in all of $(U' \times {}^*{U'})$
and thus holomorphic in $(U' \times {}^*{U'})$.
(see \cite{whitney:cav} for more).
Thus we have a real-analytic function, say
$\hat{\rho} \colon U' \to \R$
that is locally outside of $\Delta$ given by 
$i^m \prod_{j,k=1}^m (e_j(z)-\overline{e_k(z)})$.

We let $\hat{H} := \{\hat{\rho} = 0\}$.  We need to now see that $(H \cap U')
\cap \hat{H}$ is open in $H$, because then $(H \cap U') \subset \hat{H}$ as
$H$ is irreducible in $U'$
and we can apply Lemma \ref{irredlemma} as we can always pick a smaller $U'$.

It is obvious
that $\Sigma_{\gamma(0)}' \cap U'$ is in both $H$ and $\hat{H}$.  The trouble is
for other $t$, as $V \cap U'$ may in fact be empty.  Because of
how we picked $\zeta$, we note that the topological
component of $(\Sigma_{\gamma(0)}')^*$ where $\zeta$ lies is connected to 0.  So
we can find a nonsingular point $\zeta'$ of $\Sigma_{\gamma(0)}'$ on this component
that is arbitrarily close to 0, and thus inside $U'$.  We can pick a
finite sequence of overlapping neighbourhoods $\{V_j\}$ from $\zeta$ to $\zeta'$
such that inside each $V_j$, $H$ is given by $\Im f_j(z) =0$ (for
some $f_j$ holomorphic in $V_j$).
We call the final neighbourhood $V'$ and assume $V' \subset U'$
and there $H$ is given by $\Im f'(z) = 0$ (for
some $f'$ holomorphic in $V'$).  It is easy
to see that the Levi foliation is given by $f_j (z) = r$ for some real
$r$, and that these sets must agree on $V_j \cap V_k$.  Thus for some
$\epsilon'' > 0$, for all $\abs{t} < \epsilon''$, we have a component of
$\Sigma_{\gamma{t}}$ passing
through $V$ which also passes thorough $V'$ which contains $\zeta'$.  But 
$V' \subset U'$ and $\hat{H}$ and $H$ both contain all
points $\{z \mid f'(z) = t, \abs{t} < \epsilon'' \}$ and that is an open set
in $H$.

Now that we know that $H$ is contained in $\hat{H}$ we can remove
$\Delta$ which is complex analytic and work only in small neighbourhoods
where $\hat{H}$ is given by $i^m \prod (e_j(z)-\overline{e_k(z)})$.
Since $e_j(z)-\overline{e_k(z)}$ is pluriharmonic, and
thus its real and imaginary parts are pluriharmonic, meaning
that we can represent them as the imaginary part of a 
holomorphic function, that is $\Im f_{jk}(z) + i \Im g_{jk}(z)$.
Thus we get locally that
\begin{equation*}
\hat{\rho}(z,\bar{z})
=
i^m
\prod_{j,k=1}^m (\Im f_{jk}(z) + i \Im g_{jk}(z)) .
\end{equation*}
If $\Im f_{jk}(z) + i \Im g_{jk}(z)$
is zero then $(\Im f_{jk}(z))(\Im g_{jk}(z))$ is also
zero.  Thus we can make yet a larger surface by looking at
the zero set of
\begin{equation*}
\dhat{\rho}(z,\bar{z})
=
\prod_{j,k=1}^m (\Im f_{jk}(z))(\Im g_{jk}(z)) .
\end{equation*}
That is just a product of the imaginary parts of holomorphic functions.
We can now take out the set where the gradient of $f_{jk}$ and $g_{jk}$
vanish, which is a complex analytic set and work in smaller neighbourhoods
outside this set.  We can take these neighbourhoods small enough such that 
each $\Im f_{jk} = 0$ and $\Im g_{jk} = 0$ defines a nonsingular, connected
hypersurface.  The singular set of $H$ must be contained
in the intersection of at least two of these surfaces (if there is more then
one left).  This intersection is 
a generic real-analytic Levi-flat submanifold of codimension 2
outside a complex analytic subvariety by Lemma \ref{interleviflat}.
\end{proof}

%The above proof and examples known to the author seem to suggest
%the following conjecture.
%
%\begin{conj}
%Let $H \subset \C^N$ be singular irreducible
%real-analytic Levi-flat hypersurface, defined in a neighbourhood
%of the origin.  Then $H_s$ is a complex analytic subvariety.
%\end{conj}

\section{Uniqueness property for holomorphic functions} \label{sec:uniq}

For a generic submanifold $M$ through the origin in $\C^N$,
we wish to investigate when there exists a meromorphic 
function near the origin which is real valued on $M$.  By composing with a
M\"obius mapping of the real line onto the unit circle we see that this is
equivalent to the existence of a meromorphic function which is unimodular on
$M$, which in turn
means that there are two relatively prime holomorphic functions $f$ and $g$
such that on $M$, $\abs{f} = \abs{g}$.  We will thus define:
\begin{defn}
$M$ has the \emph{modulus uniqueness property} if
$\abs{f} = \abs{g}$ on $M$, for holomorphic $f$ and $g$ defined in a
neighbourhood of $M$,
implies $f = cg$ for a unimodular constant $c$.  We will say
that $M$ has the modulus uniqueness property at $p \in M$, if 
$M \cap U$ has the modulus uniqueness property for every connected
neighbourhood $U$ of $p \in M$.
\end{defn}
In the following we will denote the local CR orbit at a point $p$
by $\Orb_p$.
The motivation for our problem is the following theorem.

\begin{thm}[see \cite{BER:craut}]
\label{nicenormcoord}
Let $M \subset U \subset \C^N$ be a
generic
real-analytic nowhere minimal
submanifold of codimension $d$.
Let $p \in M$ be such that $\Orb_p$ is of maximal dimension.
Then there are coordinates
$(z,w^\prime,
w^{\prime\prime})\in\C^n\times\C^{d-q}\times\C^q=\C^N$,
where $q$ denotes the codimension of $\Orb_p$
in $M$, vanishing at $p$ such that near $p$, $M$ is
defined by
\begin{align*}
&\Im w^\prime=\varphi(z,\bar z,\Re
w^\prime,\Re w^{\prime\prime})
\\&\Im w^{\prime\prime}=0,
\end{align*}
where $\varphi$ is a real valued real-analytic function with
$\varphi(z,0,s^\prime, s^{\prime\prime})\equiv0$.
Moreover,
the local CR orbit of the point
$(z,w^\prime,w^{\prime\prime})=(0,0,s'')$, for
$s'' \in \R^q$, is given by 
\begin{align*}
&\Im w^\prime=\varphi(z,\bar z,\Re
w^\prime,s'')
\\&w^{\prime\prime}=s'' .
\end{align*}
\end{thm}

So a natural question is to ask what happens at points where $\Orb_p$ is not
of maximal dimension.  In general there do not exist local normal coordinates
such that $\Im w'' = 0$ is one of the equations for $M$,
but it is natural to ask when can we get a
meromorphic function $f$ such that $\Im f = 0$ on $M$.

Before looking at this case we summarize the results for the easy cases.
\begin{prop}
Let $M$ be a connected real-analytic CR submanifold through the origin.  Then
$M$ does not have the modulus uniqueness property at the origin
if any of the following holds,
\begin{enumerate}[(i)]
\item $M$ is not generic,
\item $M$ is totally real,
\item $M$ is nowhere minimal and $\Orb_0$ has the maximal dimension.
\end{enumerate}
On the other hand $M$ has the modulus uniqueness property at any point $p \in
M$ if
\begin{enumerate}[(i)]
\item[(iv)] $M$ is generic and minimal at some point.
\end{enumerate}
\end{prop}

\begin{proof}
The first three cases are clear.  For the last one we just note that 
if $M$ is minimal at some point, it is minimal on a dense open subset.
If we had a nonconstant meromorphic function real valued on $M$, then
on some small neighbourhood we would have that $M$ is minimal and there would
exist a holomorphic function with nonvanishing gradient which was real
valued on $M$ and this would give local foliation of $M$ by smaller
submanifolds of same CR dimension and this would violate minimality.
\end{proof}

We also note that if $M$ is not generic, but it is minimal, then 
$M$ has the modulus uniqueness property inside the intrinsic complexification
of $M$.  So since $\Orb_p$ is always minimal then if we call $\sX_p$
the intrinsic complexification of $\Orb_p$, then any meromorphic function
real valued on $M$ is constant in $\sX_p$ for any CR manifold.

It is then clearly useful to be able to construct $\sX_p$ and study its
properties.
The following constructions are described in \cite{BER:book}.
We will look at a generic submanifold $M$ defined
in normal coordinates $(z,w)$ in some 
neighbourhood $U$ of the origin, and we will assume that $U$ is small enough
such that the defining equations for $M$ complexify into $U \times {}^*U$,
and we can take $U$ to be connected.
If $M = \{ z \in U \mid r(z,\bar{z}) = 0 \}$
we let 
$\sM := \{ (z,\zeta) \in U \times {}^* U \mid r(z,\zeta) = 0 \}$.
We define the \emph{Segre manifolds} for $p \in U$
\begin{align*}
\fS_{2j+1}(p,U) :=
\{ & (z,\zeta^1,z^1,\ldots,\zeta^j,z^j) \in U \times {}^*U \times U \times
\ldots \times {}^*U\times U \mid \\
 & (z,\zeta^1),(z^1,\zeta^1),\ldots,(z^j,\zeta^j),
 (z^1,\zeta^2),\ldots,(z^{j-1},\zeta^j),(z^j,\bar{p}) \in \sM \}
\end{align*}
and
\begin{align*}
\fS_{2j}(p,U) :=
\{ & (z,\zeta^1,z^1,\ldots,z^{j-1},\zeta^j) \in U \times {}^*U \times U \times
\ldots \times U \times {}^*U \mid \\
 & (z,\zeta^1),(z^1,\zeta^1),\ldots,(z^{j-1},\zeta^{j-1}),
 (z^1,\zeta^2),\ldots,(z^{j-1},\zeta^j),(p,\zeta^j) \in \sM \}
\end{align*}
where $\fS_1(p,U) = \{ z \in U \mid (z,\bar{p}) \in \sM \}$.
If we define $\pi \colon \C^N \times \ldots \times \C^N \to \C^N$ be the
projection to the first coordinate, then we can define the
\emph{Segre sets} for $p \in U$ by
$S_k(p,U) := \pi (\fS_k (p,U) )$.  Note that both $S_k(p,U)$
and $\fS_k(p,U)$ depend on both the point $p$ and the neighbourhood $U$.

We have the following proposition, of which the
first part is proved in \cite{BER:book}
(Proposition 10.2.7), second part is then immediate.

\begin{prop} \label{inductivedefS}
For $k \geq 1$ we have
\begin{equation*}
S_k(p,U) =
\bigcup_{q \in S_{k-1}(p,U)}
S_1(q,U) ,
\end{equation*}
and if $k \geq 2$ we have
\begin{equation*}
S_k(p,U) =
\bigcup_{q \in S_{k-2}(p,U)}
S_2(q,U) .
\end{equation*}
\end{prop}

Further, for normal coordinates where $U = U_z \times U_w$
we have the following (again proved in \cite{BER:book} as part of
Proposition 10.4.1):

\begin{prop} \label{S2normcoord}
Let $M$ be given by $w = Q(z,\bar{z},\bar{w})$
in normal coordinates in $U$ and let $p = (z^0,w^0)$.
Then there exists
an open set $V \subset {}^* U_z$ ($0 \in V$) such
that $(z,w) \in U$ is in $S_2(p,U)$ if
and only if there exists $\zeta \in V$
such that $w = Q(z,\zeta,\bar{Q}(\zeta,z^0,w^0))$.
\end{prop}

The above $V$ is the set of all $\zeta \in {}^* U_z$
such that $\bar{Q}(\zeta,z^0,w^0) \in {}^*U_w$.  In particular $0 \in V$.
With this we prove the following useful lemma.

\begin{lemma} \label{S2noz}
Suppose that $M \subset U \subset \C^N$
is a generic submanifold
given by normal coordinates defined near
the origin for a suitable $U$.
Then for any point $p = (z^0,w^0) \in U$, the variety
$\{ (z,w) \in U \mid w = w^0 \}$ is contained inside $S_2(p,U)$ (the second
Segre set at $p$).
\end{lemma}

\begin{proof}
Let $M$ be given by $\{ w = Q(z,\bar{z},\bar{w}) \}$ in normal coordinates.
Thus
$S_2(p,U) = \{ (z,w) \mid w = Q(z,\zeta,\bar{Q}(\zeta,\bar{z}^0,\bar{w}^0)),
\zeta \in V \}$, where $V$ is as in Lemma \ref{S2normcoord}.
In particular $0 \in V$ and thus since we are in normal
coordinates, $Q(z,0,w) \equiv Q(0,z,w) \equiv w$.  Thus 
$\{ (z,w) \mid w = w^0 \} \subset S_2(p,U)$.
\end{proof}

To be able to use this we note the following theorem given and proved in
\cite{BER:book} (Theorems 10.5.2 and 10.5.4).

\begin{thm} \label{SpXp}
If $M$ is as above then there exists a number $j_0$ such that
for every sufficiently small neighbourhood $U$ of $p \in M$,
$S_{2j_0}(p,U)$ coincides with $\sX_p$ as germs at $p$,
the intrinsic complexification of $\Orb_p$.
\end{thm}

The number $j_0$ is called the Segre number of $M$ at $p$,
but we are only interested in the fact that such a number exists and not how
it is arrived at.  Another useful proposition from \cite{BER:book}
(Proposition 10.2.28)
is the following.

\begin{prop} \label{Sconn}
Let $p \in M \subset U$ and an integer $k_0 \geq 1$.
Then there exist neighbourhoods $U'' \subset U'
\subset U$ of $p$ such that for all $q \in U''$, $\fS_k(q,U')$ is connected
for all $k \leq k_0$.
\end{prop}

Next we assume that $U$ is sufficiently nice (for example a polydisc).

\begin{lemma} \label{Xpnoz}
Given $M \subset U$ in normal coordinates, then there is a small neighbourhood
of the origin $V$ such that for $p \in M \cap V$,
$\sX_p$ contains $\{ (z,w) \in U \mid w = w^0 \}$
as germs at any $(z^0,w^0) \in \sX_p$.
If $\sZ_p$ is the intersection of all complex subvarieties of $U$ which contain
$\sX_p$, then $\sZ_p$ contains $\{ (z,w) \in U \mid w = w^0 \}$
for any $(z^0,w^0) \in \sZ_p$.
\end{lemma}

\begin{proof}
Let $M$ be in normal coordinates.
We can always take $U$ to be even smaller, so
by Proposition \ref{Sconn} for a small
enough neighbourhood of the origin $U$,
there is a yet smaller neighbourhood of the origin $V$ such that
for $p \in V$, $\fS_k(p,U)$ is connected, for $k \leq 2(d+1) + 2$,
$d$ being the codimension of $M$.
Note that the Segre number of $M$ at any point is always less then or
equal to $d+1$.
By Theorem \ref{SpXp} we know
$S_{2(d+1)}(p,W) = \sX_p$ as germs for some small neighbourhood $W$ of $p$.
Hence $S_{2(d+1)+2}(p,W) = \sX_p = S_{2(d+1)}(p,W)$ as germs at $p$.
Let $k = 2(d+1)$.
By Proposition \ref{inductivedefS}, 
$S_{k+2}(p,U)$ is a union of $S_2(q,U)$ for $q \in S_{k}(p,U)$,
and by Lemma \ref{S2noz} each $S_2(q,U)$ contains the set $\{ (z,w) \mid w =
w(q) \}$.
In particular $S_{k+2}(p,U)$ contains the set
$\{ (z,w) \mid w = w(q) \}$ for each $q \in \sX_p$ (for some small enough
representative of the germ $\sX_p$).
Now we note that
$\fS_{k+2}(p,W)$ is an open submanifold of $\fS_{k+2}(p,U)$,
which is connected.  We pull back the mapping $(z,w) \mapsto z$ to
$\fS_{k+2}(p,U)$ and look at its rank to conclude that
for $q \in S_{k+2}(p,W)$ we have
$\{ (z,w) \mid w = w(q) \} \subset S_{k+2}(p,W)$ as germs at $q$.  This
proves the first part.

%So the first part follows by
%putting Proposition \ref{inductivedefS}, Lemma \ref{S2noz}
%and Theorem \ref{SpXp} together.
To see the second part suppose that $\sZ_p$ did depend on $z$.
Then we can intersect $\sZ_p$ with $\{ (z,w) \mid z = z^0 \}$
and the intersection must still contain $\sX_p$ projected on the $w$ coordinate
(it is of the form $\sX_p = \C_z \times (\sX_p)_w$).
So we would get a different
complex variety $\sZ_p'$ which contains $\sX_p$.  Intersection 
of $\sZ_p$ and $\sZ_p'$
would violate minimality of $\sZ_p$.
\end{proof}

\begin{thm} \label{meronoz}
Suppose that $M$ is generic in normal coordinates.  Suppose that $f$ and
$g$ are two holomorphic functions such that $\abs{f} = \abs{g}$ on $M$.
Then $f/g$ depends only on $w$;
in other words, if $h$ is a meromorphic function
which is real valued on $M$, then $h$ depends only on $w$.
\end{thm}

\begin{proof}
Obviously we only need to prove the first part as the second part follows.
We can work in arbitrarily small neighbourhood $U$ of the origin.
As we noted before since $\Orb_p$ is minimal in $\sX_p$ we know that $f = cg$
in $\sX_p$ for any point $p$ (where $c$ depends on $p$ of course).
That is that the function $f/g$ is constant on $\sX_p$ (if we take $p$ outside
the zero set of $g$).
Since we know that as germs 
$\{ (z,w) \mid w = w^0 \} \subset \sX_p$, then for any $1 \leq j \leq n$
we have
$\frac{\partial}{\partial z_j} (f/g) = 0$ at $p$.  Since $M$ is generic and
since $g=0$ is a proper subvariety of $M$, then
$\frac{\partial}{\partial z_j} (f/g) = 0$ 
holds
for an open set of $p$ in $M$, and then it holds for an open subset of $U$
and thus for all of $U$.
\end{proof}

\section{Submanifolds inside Levi-flat hypersurfaces} \label{sec:flatness}

Since the question of the modulus uniqueness property of $M$
(or alternatively of existence of a meromorphic function which is real
valued on $M$)
is the same as a question of $M$ being contained in a certain kind of possibly
singular real-analytic Levi-flat hypersurface, we can ask a weaker question;
when is $M$ contained in any possibly singular real-analytic Levi-flat
hypersurface?  We will consider $M$ to be inside a hypersurface $H$
if $M \subset \overline{H^*}$.

%(FIXME: I thought this could be generalized to lower codimension,
% it can a little if we assume M intersects H^*)
\begin{prop}
Suppose $M$ is a connected generic real-analytic submanifold of
codimension 2 in normal coordinates $(z,w)$
and $M \subset \overline{H^*}$ where $H$ is a
irreducible possibly singular real-analytic
Levi-flat hypersurface.
Then in a possibly smaller neighbourhood of the origin,
there exists a Levi-flat hypersurface $\hat{H}$ defined by
$\{ (z,w) \mid \rho(w,\bar{w}) = 0 \}$ such that
$M \subset \overline{\hat{H}^*}$ as germs at 0.
Furthermore, if $M$ is not Levi-flat then $H = \hat{H}$ as germs at 0.
\end{prop}

\begin{proof}
If $\Orb_p$ is constantly of codimension 2 in $M$
or constantly of codimension 1 in $M$, then by Theorem \ref{nicenormcoord}
we have a holomorphic function near the origin which is real valued on $M$
and thus by Theorem \ref{meronoz} the defining equation for $M$ already
does not depend on $z$.

By Corollary \ref{minhnotminimal}, $M$ cannot be minimal at any point.
So suppose that $M$ is not minimal and $\Orb_p$ is not of constant dimension.
This means that it is not Levi-flat and thus by Corollary
\ref{flathypersingcor} it cannot be contained in $H_s \cap \overline{H^*}$
and thus must intersect $H^*$.  This means that it must in fact intersect
$H^*$ on a dense open set in $M$ (as $H_s$ is contained in a
proper subvariety).
Suppose $H$ is defined in $U$
by $\{ \rho(z,w,\bar{z},\bar{w}) = 0 \}$, in particular $H$ is closed in $U$.
Then for $p \in M \cap H^*$ we can see that $\sX_p \subset H$,
since in small enough neighbourhood of $p$, such as
we have by Theorem \ref{SpXp}, the $k$th Segre set
of $M$ is contained in the $k$th Segre set of $H$, and the Segre sets
of $H$ all lie in $H$ for small enough neighbourhood of a nonsingular point
of $H$.  By Lemma \ref{sprimenonempty}, the Segre variety of $H$
at $p$ agrees with the Levi foliation of $H$ at $p$, and since this (the
Segre variety of $H$) is a proper
subvariety of $U$, then if $\sZ_p$ is the smallest complex subvariety
of $U$ which contains $\sX_p$, then $\sZ_p \subset H$.
This means in particular that
$\big( \C_z \times pr_w(M \cap H^*) \big) \cap U \subset H$ (where $pr_w$ is
the projection onto second factor in the normal coordinates $(z,w)$), since
$\sZ_p$ contains
all the $(z,w) \in U$ for fixed $w$ by Lemma \ref{Xpnoz}.
As $H$ is closed and
$M \cap H^*$ is dense in $M$, then
$\C_z \times pr_w(M) \subset H$.  Fix $z^0$ such that
$\rho(z^0,w,\bar{z}^0,\bar{w}) = 0$ defines a hypersurface in $\C_w$,
then this hypersurface is Levi-flat in $\C_w$.
Define $\hat{H}$ by 
$\{ (z,w) \mid \rho(z^0,w,\bar{z}^0,\bar{w}) = 0 \}$, this is 
Levi-flat again and further $M \subset \hat{H}$.

It is then clear that since $\sX_p \subset H$, then $\sX_p \subset \hat{H}$,
thus near points $p$ where $\Orb_p$ is of codimension 1 in $M$, they locally
give a branch of a nonsingular Levi-flat hypersurface which must be contained in
$\hat{H}$.  Thus $M \subset \overline{\hat{H}^*}$.

If $M$ is not Levi-flat then uniqueness of $H$ comes from
the fact that if $M$ would be contained in two different Levi-flat
hypersurfaces say $H$ and $H'$
it would be contained in their intersection and
thus would be contained in the singular set of $H \cup H'$
and this is impossible by Corollary \ref{flathypersingcor}.
\end{proof}

Our method of looking at projections onto the second factor of
normal coordinates yields
also first part of Theorem \ref{fullprojwflatthm} which we can state
as follows.

\begin{thm} \label{projwflatthm}
Let $M$ be a germ of a
generic real-analytic
codimension 2 submanifold through the origin
given in normal coordinates $(z,w)$ and let $M$ be nowhere minimal.  
Then $M \subset \overline{H^*}$, where $H$ is a germ of
a possibly singular real-analytic
Levi-flat hypersurface if and only if
the projection of $M$ onto the second factor in $(z,w)$
is contained in a germ of a
possibly singular real-analytic hypersurface.
Moreover, if $M$ is not Levi-flat, then $H$ is unique.
\end{thm}

Note that this theorem also gives a test for certain submanifolds being
nowhere minimal.  If we can compute a hypersurface containing the
projection of $M$ to the $w$ coordinate, we need only check if it is
Levi-flat or not.

\begin{proof}%[Proof of Theorem \ref{projwflatthm}]
The forward direction and uniqueness
is proved by the preceding proposition.  So suppose
that $M \subset H$ where $H = \C_z \times H_w$ is a possibly singular
hypersurface.  We can assume that $H$ is irreducible.

First suppose that $\Orb_0$ is of maximal
dimension.  Then by Theorem \ref{nicenormcoord} there exists (near 0)
a holomorphic
function real valued on $M$ which thus defines a Levi-flat hypersurface
(nonsingular one in fact).  Also by Theorem \ref{meronoz} this function
only depends on the $w$ coordinate, this means that it really defines a
Levi-flat hypersurface in $\C^d$ (the $w$ space)
that contains $pr_w(M)$.

Next suppose that $\Orb_0$ is not of maximal dimension.
Fix a certain neighbourhood $U$ where $M$ is defined in the given
normal coordinates.  By Proposition \ref{Sconn} we can then pick a
smaller $0 \in U' \subset U$ such that for all $p \in U'$, the Segre manifold
$\fS_k(q,U)$ is connected.
Making $U'$ smaller we can assume it is
of the form $U'_z \times U'_w$ where both $U'_z$ and $U'_w$ are polydiscs.
We will pick a point $p \in M \cap U'$ where $\Orb_p$ is of maximal dimension
(of codimension 1 in $M$).

By choosing above $U$ small enough we can ensure that $pr_w(M)$
is subanalytic (see \cite{BM:semisub}).  We look at a nonsingular point
of this projection of highest dimension in $pr_w(M) \cap U'_w$.
Obviously at this point $pr_w(M) \cap U'_w$ is either a hypersurface or
codimension 2 since it is contained in $H_w$.  If $pr_w(M)$ was a codimension 2
submanifold near some point, then it would be totally real, and thus $M$ above
it would be Levi-flat which is not the case.  Thus there must be nonsingular
points where $pr_w(M) \cap U'_w$ is a hypersurface.
Further, since the $\sX_q$ really
depends only on the $w$ variables,
it is clear that there is a point $p \in M \cap U'$,
such that $pr_w(p)$ is a nonsingular point of $pr_w(M) \cap U'_w$,
and such that $\Orb_p$ is
of maximal
dimension.  Next, pick a small enough neighbourhood $V \subset U'$ of $p$,
such that $pr_w(M \cap V)$ is a nonsingular hypersurface.  Then 
$pr_w(M \cap V)$ agrees with one of the branches of $H_w$ at $pr_w(p)$.

Locally in $V$ (possibly taking smaller $V$)
again we have a holomorphic function $f$ in a neighbourhood of $p$
that is real valued on $M$.  We notice that in the proof of Lemma \ref{Xpnoz}
the only reason why we restrict to a smaller neighbourhood is so that we
can apply Proposition \ref{Sconn}, and hence we could have picked a
neighbourhood of any point in $U$.  So we see that in the proof
of Theorem \ref{meronoz} we did not need to pick a neighbourhood of the
origin, but we could have just used $V$ as given above (possibly making
it smaller).
Hence
$f$ only depends on $w$, and thus again $\Im f = 0$ defines
a Levi-flat hypersurface near $p$ which contains $M$ near $p$.  So in
$\C^d$ (the $w$ coordinates)
this hypersurface contains $pr_w(M \cap V)$ and thus agrees with
a branch of $H_w$ near $pr_w(p)$.  By Lemma
\ref{flatsomewhereflateverywhere}, $H_w$ must be a Levi-flat
hypersurface, and we are done.
\end{proof}

\section{Almost minimal submanifolds} \label{sec:almostmin}

As we have already seen, if $M \subset \overline{H^*}$ and $M$ is a generic
nowhere minimal
codimension 2 real-analytic submanifold,
then at a point $p \in M \cap H^*$ where $\Orb_p$ is of codimension
1 in $M$, $\sX_p \subset H^*$, that
is, $\sX_p$
gives the Levi foliation of $H$.  By Lemma \ref{sprimenonempty}, we have
that locally the Segre variety $\Sigma_p$
of $H$ in $U$ contains $\sX_p$, and for $p \in
H^*$, $\Sigma_p$ is a proper analytic subvariety of $U$.  So an obvious condition
for $M$ to be contained in a Levi-flat hypersurface is that $\sX_p$ is
contained in a proper complex subvariety of $U$.
Since $\sX_p$ is the smallest germ of a complex variety
containing $\Orb_p$, we let
$\sZ_p = \sZ_{U,p}$ be the smallest complex subvariety of $U$
that contains $\Orb_p$ (and thus $\sX_p$).

\begin{defn}
Let $M \subset U \subset \C^N$ be a generic submanifold.  We will say
that $M$ is \emph{almost minimal in $U$}, if there exists a point $p$
such that $\sZ_{U,p}$ contains an open set, and we will say that $p$ makes
$M$ almost minimal in $U$.
We'll say that a generic submanifold
$M$ is \emph{almost minimal at $p$}, if it is almost minimal
in every neighbourhood of $p$.
\end{defn}

If $M$ is minimal at $p \in U$, then it is, of course, almost
minimal in $U$.  And if a connected $M$ is real-analytic and minimal at
one point, it is minimal on an open dense set, and thus it is almost minimal
at every point.

An example of a nowhere minimal submanifold that is almost minimal
is the $M_\lambda$ family given in the introduction for $\lambda$ irrational.
See \S \ref{sec:ex} for this example worked out.  It should be noted that if
$M$ is nowhere minimal, then the points where it is almost minimal are
contained in a proper real analytic subvariety in $M$.
This is because if $M$ is almost minimal at $p$
and nowhere minimal, then $\Orb_p$ must not be of maximal dimension.

\begin{thm} \label{m4flatthm}
Suppose that $M \subset \C^N$ is a germ of a
real-analytic generic submanifold of
codimension 2 through 0,
and suppose $M \subset \overline{H^*}$
where $H$ is a germ of a possibly singular real-analytic Levi-flat
hypersurface.
Then $M$ is not almost minimal at 0.
\end{thm}

\begin{proof}
Let $U$ be a small enough connected neighbourhood of the origin such that
both $M$ and $H$ are closed in $U$ and further such that their defining
equations are complexifiable in $U$.
$M$ cannot be minimal at any point by Theorem \ref{minhnotminimal}.
Further, if $M$ is Levi-flat then
$\Orb_p$ is constantly of codimension 2 in $M$.  This means that $\Orb_p$
is in fact complex analytic and is contained in the Segre variety (the first
Segre set of $M$ in $U$) and thus cannot be almost minimal.

So suppose on a dense open set of points of $M$, $\Orb_p$ is of codimension
1 in $M$, and in fact, if $p$ makes $M$ almost minimal in $U$
then $\Orb_p$ has to be of codimension 1 in $M$.
Further, $M \cap H^*$ is non-empty
(since $M$ is not Levi-flat)
and as noted before is thus open and dense in $M$.  Also as noted above, the
$p$ that makes $M$ almost minimal cannot lie in $M \cap H^*$.

So pick a small neighbourhood of any $p \in M \cap H_s$
where $\Orb_p$ is of codimension 1 in $M$.  Then by
Theorem \ref{nicenormcoord}, there is a small neighbourhood
$V$ of $p$ on which there exist normal coordinates
$(z,w) \in \C^n \times \C^2$
vanishing at $p$, such
that $M$ is given by
$\Im w_1 = \rho(z,\bar{z},\Re w)$ and $\Im w_2 = 0$,
and further, that the $\sX_q$ are then given by
$w_2 = s$ (we'll denote this set as $\{w_2 = s \}$)
for some $s \in (-\epsilon,\epsilon)$.
We can take $V$ to be a polydisc in the $(z,w)$ coordinates.
If $M \cap \{ w_2 = s \}$
(which is the CR orbit) contains a point which is in $H^*$,
then as we reasoned above $\{ w_2 = s \} \subset H$ since it agrees
with the Levi foliation of $H$ at some point in $H^*$.  
As $M \cap H^*$ is dense in $M$, then $\{ w_2 = s \} \subset H$
for all $s \in (-\epsilon,\epsilon)$.  This means that in
$V$, $\Im w_2$ divides the defining function of $H$ in $U$.
Thus the Segre variety of $H$ in $U$ contains the Segre variety
of $\{ \Im w_2 = 0 \}$ at all points in
$\{ \Im w_2 = 0 \}$.  We wish to show that $\Orb_p$ is contained
in a proper complex analytic subvariety.  Either it is contained
in a nondegenerate Segre subvariety of $H$ in $U$ or the Segre variety of $H$
in $U$ is degenerate at all points of $\Orb_p = M \cap \{ w_2 = 0 \}$,
but the set of
points where the Segre variety of $H$ is degenerate is a proper analytic
subset as we remarked before.
In any case $p$ does not make
$M$ almost minimal in $U$, and thus $M$ is not almost minimal in $U$.
\end{proof}

\begin{cor}
Suppose that $M \subset \C^N$ is a connected
real-analytic generic submanifold of
codimension 2 through 0, and $M$ is almost minimal at $0$.
Then $M$ has the modulus uniqueness property at $0$.
\end{cor}

\section{Infinitesimal CR automorphisms} \label{sec:craut}

We will now look at the dimension of $\hol(M,p)$, the space of infinitesimal
holomorphisms at $p$ (see the introduction for terminology) if $M$ is almost
minimal at $p$.  As motivation 
we have the following theorem.

\begin{thm}[Baouendi-Ebenfelt-Rothschild see \cite{BER:craut}]
\label{berholfindim}
Let $M \subset \C^N$ be a connected real-analytic CR submanifold that is
holomorphically nondegenerate.
If $M$ is minimal at any point $p \in M$, then $\dim_\R \hol(M,q) < \infty$ for
all $q \in M$.  If $M$ is nowhere minimal then
$\dim_\R \hol(M,q) = 0$ or
$\dim_\R \hol(M,q) = \infty$
for $q$ in a dense open subset of $M$.
\end{thm}

Thus it remains to see at exactly what points $\hol(M,q)$ is finite dimensional
in case $M$ is nowhere minimal.  Our main result of this section is that
it turns out that the points where $M$ is almost minimal are such points.
We restate Theorem \ref{holfindim} from the introduction for convenience.

\begin{thmnonum}
%\label{holfindim}
Let $M \subset \C^N$ be a connected, real-analytic holomorphically
nondegenerate generic submanifold and suppose
$p \in M$ and $M$ is almost minimal at $p$.  Then
\begin{equation*}
\dim_\R \hol(M,p) < \infty .
\end{equation*}
\end{thmnonum}

The proof is essentially the same as in \cite{BER:craut} or \cite{BER:book}
for minimal submanifolds,
although we will require Lemma \ref{keylemma} to modify that proof.
It would not be needed, if we had a more general way of showing that certain CR orbits
(of the highest dimension for example) were holomorphically nondegenerate
whenever $M$ was.  However this is not so.  For example, the manifold
defined in $(z_1,z_2,w_1,w_2) \in \C^4$, by 
\begin{align*}
& \Im w_1 = \abs{z_1}^2 + (\Re w_2) \abs{z_2}^2 ,
\\
& \Im w_2 = 0 ,
\end{align*}
is holomorphically nondegenerate.  However, the CR orbit at 0 is defined
by $\Im w_1 = \abs{z_1}^2$ and $w_2 = 0$, and so $\frac{\partial}{\partial
z_2}$ is a holomorphic vector field tangent to it.  We can, however,
prove
the following result for almost minimal submanifolds.

\begin{lemma}
Suppose $M \subset U$ is a holomorphically nondegenerate generic submanifold,
and $p \in M$ is such that $\sZ_{U,p} = U$, that is, $p$ makes $M$ almost
minimal in $U$.  Then $\Orb_p$ is holomorphically
nondegenerate.
\end{lemma}

The proof is essentially contained the proof of Theorem
\ref{holfindim} below, and uses the following technical lemma.

\begin{lemma}
\label{keylemma}
Let $M \subset U_z \times U_w \subset \C^N$ be a generic submanifold given
in normal coordinates $(z,w)$ in a sufficiently small $U = U_z \times U_w$
by $w = Q(z,\bar{z},\bar{w})$. 
Suppose there exists a holomorphic function $f \colon \C^n \times \C^n
\times \C^d \to \C$ defined in a neighbourhood of the origin such that
$f(z,\bar{z},w)$ is defined in $U$, and there exists a point $p \in
M$ and $f(z,\bar{z},w) = 0$ on $\Orb_p$.  Then there exists a
holomorphic function $g \colon U_w \subset \C^d \to \C$
such that $g(w) = 0$ on $\Orb_p$.
\end{lemma}

\begin{proof}
First note that Lemma \ref{Xpnoz} implies that locally near $p$, we can find
a germ of a holomorphic function $\varphi$ such that $\varphi(w) = 0$ defines
$\sX_p$.  Thus we can do a local change of coordinates in $w$ only,
setting $w' = \psi(w)$ and $w'' = \varphi(w)$ for some function $\psi$.
So locally we have $\Orb_p$ defined in the coordinates $z,w',w''$ (which are no
longer normal coordinates) by $w' = \tilde{Q}(z,\bar{z},\bar{w}')$
and $w'' = 0$, for some function $\tilde{Q}$ defined in a neighbourhood
of $p$.
We can now also write $f$ in the $z$, $w'$, $w''$ coordinates by abuse of
notation as
$f(z,\bar{z},w',w'')$.
Assuming $U$ is small enough and the neighbourhood where $w',w''$ are defined
is also small enough we can
define a complexified version of $\Orb_p$ by setting $\bar{w}' = \xi$
and $\bar{z} = \zeta$ by $\xi = \bar{\tilde{Q}}(\zeta,z,w')$ and call this
$\sC$.  Since $f(z,\bar{z},w',0) = 0$ on $\Orb_p$, then as
$\Orb_p$ is maximally real in $\sC$ we have that $f(z,\zeta,w',0) = 0$ on
$\sC$ and as $z,\zeta,w'$ are free variables on $\sC$ we know that 
$f(z,\bar{z},w',w'') = 0$ when $w'' = 0$, but $w'' = 0$ defines $\sX_p$,
so $f$ is identically zero on all of $\sX_p$.  Since $\sX_p$ is defined
by an equation which is independent of $z$, then if we 
fix $z^0$ where $(z^0,w^0) = p \in M$, and we take
$g(w) := f(z^0,\bar{z^0},w)$, then $g(w)$ as a function of
$(z,w)$ but independent of $z$ is zero on $\sX_p$ and thus on $\Orb_p$.  And
$g$ is defined in all of $U_w$ and thus we are done.
\end{proof}

We need to characterize $\hol(M,p)$ in a more natural way for the proof and
the following proposition is proved in \cite{BER:book} (Proposition
12.4.22).

\begin{prop}
\label{holcrautchar}
Let $M \subset \C^N$ be a real-analytic generic submanifold and $X$ a germ of
a real, real-analytic vector field on $M$.  Then $X \in \hol(M,p)$ if and only
if there exists a germ $\sX$ at $p$ of a holomorphic vector field in $\C^N$
such that $\Re \sX$ is tangent to $M$ and $X = \Re \sX|_M$.
\end{prop}

It is not hard to see that
if $\sX$ is a holomorphic vector field as above and $\tilde{\varphi}(z,\tau)$
is the holomorphic flow of $\sX$ and $X = \Re \sX$ and
$\varphi(z,\bar{z},t)$ is the flow of $X$, then $\varphi$ and
$\tilde{\varphi}$ coincide when $t = \tau \in \R$.

We will need the notion of $k$-nondegeneracy, but instead
of giving the definition of being $k$-nondegenerate at a point, we
can just take the following proposition from \cite{BER:book}
(Corollary 11.2.14) and treat it as a definition.

\begin{prop}
\label{knongendef}
Let $M \subset \C^N$ be a real-analytic generic submanifold of codimension
$d$ and CR dimension $n$ given in normal coordinates $Z = (z,w) \subset U
\subset \C^n\times\C^d$ by
$w = Q(z,\bar{z},\bar{w})$.  Then $M$ is $k$-degenerate at $p = (z_p,w_p)$
(sufficiently close to 0)
if and only if
\begin{equation*}
\operatorname{span}
\left\{
\left(\frac{\partial}{\partial \bar{z}} \right)^\alpha
\frac{\partial \bar{Q}_j}{\partial Z} (\bar{z}_p,z_p,w_p)
\mid j = 1,\dots,d, 0 \leq \abs{\alpha} \leq k
\right\} = \C^N .
\end{equation*}
\end{prop}

We must prove a result about finite jet determination of
biholomorphisms of almost minimal submanifolds, which may be of interest on
its own.
As before let $\sZ_{U,p}$ be the smallest complex analytic variety containing
$\Orb_p$.  So if $M$ is almost minimal in $U$ and $p$ is the point that
makes it almost minimal then we have the following proposition.

\begin{prop}
\label{finitedetthm}
Let $M,M' \subset \C^N$ be real-analytic generic submanifolds of codimension
$d$ defined
in open sets $U$ and $U'$ respectively.  Let $f$ and $g$ be two
holomorphic mappings taking $U$ to $U'$ and $M$ to $M'$.
Let $p \in M$ be such that
$\sZ_{U,p} = U$ and suppose $M$ is
$k_0$-nondegenerate at $p$.  Also suppose that $f(p) = g(p) = p'$,
$f_*(T^c_{p}M) = T^c_{p'}M$ and
$g_*(T^c_{p}M) = T^c_{p'}M$.  Then if
$j^{(d+1)k_0}_{p} f = j^{(d+1)k_0}_{p} g$ then $f = g$.
\end{prop}

\begin{proof}
The proposition follows from Corollary 12.3.8 in \cite{BER:book} which is
a slightly stronger result than the above, but which says
that $f = g$ in $\Orb_p$ only.  As $\sZ_{U,p} = U$, then
of course $f = g$ everywhere on $U$.
\end{proof}

To be able to use Proposition \ref{finitedetthm}
we need to know that $M$ is $k$-nondegenerate at the right points.  From
\cite{BER:book} we have the following lemma (part of Theorem 11.5.1).

\begin{lemma}
\label{holnongenknongen}
Suppose $M \subset \C^N$ is a connected real-analytic generic submanifold
of CR dimension $n$ 
that is holomorphically nondegenerate.
Then there exists a proper real-analytic subvariety $V \subset M$
such that
$M$ is $\ell$-nondegenerate for all $p \in M \backslash V$ for some
$1 \leq \ell \leq n$.
\end{lemma}

The $\ell = \ell(M)$ is the \emph{Levi-number} of $M$.

\begin{proof}[Proof of Theorem \ref{holfindim}]
First suppose that $X^1,\ldots,X^m \in \hol(M,p)$ are linearly independent
over $\R$.  Suppose that $x = (x_1,\ldots,x_r)$ be local coordinates for $M$
vanishing at $p$.  Here we may write $X^j = \sum_{k=1}^{r}
X^j_k(x)\frac{\partial}{\partial x_k}$, or for short
$X^j \cdot \frac{\partial}{\partial x}$.  We let $y \in \R^m$ and denote
by $\varphi(t,x,y)$ the flow of the vector field $y_1X^1 + \cdots + y_mX^m$,
that is the solution of
\begin{align*}
& \frac{\partial \varphi}{\partial t} (t,x,y) =
\sum_{j=1}^m y_j X^j(\varphi(t,x,y)) ,
\\
& \varphi(0,x,y) = x .
\end{align*}

Since $\varphi(st,x,y) = \varphi(t,x,sy)$ (which follows from
the uniqueness of the solution), we can choose $\delta > 0$
small enough such that there exists $c > 0$ such that
the flow is smooth for $(t,x,y)$ where $\abs{t} \leq 2$,
$\abs{x} \leq c$ and $\abs{y} \leq \delta$.  We look at
the time-one map denoted by
\begin{equation*}
F(x,y) := \varphi(1,x,y) .
\end{equation*}
We have the following
lemma proved in \cite{BER:craut} and \cite{BER:book} (Lemma 12.5.10).

\begin{lemma}
\label{flowlemma}
Let $F$, $x$, $c$, and $\delta$ be as above.
There exists $\gamma > 0$ such that $\gamma < \delta$ such that
for any fixed $y^1, y^2 \in \R^m$ where $\abs{y^j} \leq \gamma$, $j =1,2$,
if $F(x,y^1) \equiv F(x,y^2)$ for $\abs{x} \leq c$ then necessarily
$y^1 = y^2$.
\end{lemma}

Suppose that $X^j$ are as above and are in $\hol(M,p)$.
Denote by $V \subset M$ the
neighbourhood of $p$ given by $\abs{x} < c$ where $x$ and $c$ are as above.
Let $\gamma > 0$ be picked as in Lemma \ref{flowlemma}.
From Proposition \ref{holcrautchar} (and discussion afterward)
it follows that for a fixed $y$ such
that $\abs{y} < \gamma$, there exists
a biholomorphism $z \mapsto \tilde{F}(z,y)$ defined in some connected open
neighbourhood $U \subset \C^N$ of $V \subset M$,
taking $M$ into $M$.  We can take $\gamma$ smaller if necessary.
Further, if $z(x)$ is the
parametrization of $M$ near $p$, these satisfy $F(x,y) = \tilde{F}(z(x),y)$,
where $F$ is the time-one map defined above.

As $M$ is holomorphically nondegenerate,
then by Lemma \ref{holnongenknongen} we have that outside a real-analytic
set it is $\ell$-nondegenerate.  Note that by Proposition \ref{knongendef}
we have that this set is actually contained in a set
defined by the vanishing of
a function of the
form $\varphi(z,\bar{z},w)$, that is a real-analytic function in $z$,
but holomorphic in $w$.  Since $M$ is almost minimal at $p$
and if $\sZ_{U,q} = U$, then
$\sZ_{U,q'} = U$ for all $q' \in \Orb_{q}$, we know by
Lemma \ref{keylemma} that there must exist
a $q \in M$ such that $\sZ_{U,q} = U$ and $M$ is $\ell$-nondegenerate
at $q$.

Thus we have satisfied requirements of Proposition
\ref{finitedetthm}, and by applying Lemma \ref{flowlemma} we see that we
have an injective
mapping
\begin{equation*}
y \mapsto j^{(d+1)\ell}_q \tilde{F} (\cdot,y) \in J^{(d+1)\ell}(\C^N,\C^N)_q ,
\end{equation*}
where $J^{(d+1)\ell}(\C^N,\C^N)_q$ is the jet space at $q$ of germs
of holomorphic mappings from $\C^N$ to $\C^N$.  As 
$J^{(d+1)\ell}(\C^N,\C^N)_q$ is finite dimensional, then 
obviously $m \leq \dim_\R J^{(d+1)\ell}(\C^N,\C^N)_q$.  Thus
$\dim_\R \hol(M,p) \leq \dim_\R J^{(d+1)\ell}(\C^N,\C^N)_q$.
\end{proof}

\section{Algebraic submanifolds} \label{sec:alg}

A manifold is real-algebraic if it is contained in a real-algebraic variety
of the same dimension.  The following theorem is basically
proved in \cite{BER:book} (Theorem 13.1.10).
It is also easily seen as a direct consequence of Tarski-Seidenberg (see
\cite{BM:semisub}) and
of the Chevalley theorem (see for example \cite{LS:book}).
That is, projections of real or complex algebraic varieties
are either semi-algebraic (in the real case)
or constructible (in the complex case) but in both cases they are contained
in a real or complex algebraic variety of the same dimension.  And since
$\sX_p$ is locally given as projection of a
Segre manifold which is complex-algebraic if $M$ is real-algebraic, we have the
following.

\begin{thm}
Let $M$ be a real-algebraic generic submanifold and $p \in M$.
Then $\Orb_p$ is real-algebraic and similarly $\sX_p$
is contained in a complex algebraic variety of the same dimension.
\end{thm}

If $M$ is nowhere minimal and real-algebraic, $U$ is
an open set, and $p \in M \cap U$,
then $\Orb_p$ is contained in a proper complex analytic
subvariety of $U$.  So we have the following corollary.

\begin{cor}
Suppose $M \subset \C^N$
is a connected real-algebraic generic submanifold.  Then
$M$ is almost minimal at $p \in M$
if and only if $M$ is minimal at some point.
\end{cor}

\begin{cor}
Let $M$ be nowhere minimal real-analytic generic submanifold which
is almost minimal at $p \in M$.  Then $M$ is not biholomorphic to a
real-algebraic generic submanifold.
\end{cor}

This is because almost minimality would be preserved under biholomorphisms.
The $M_\lambda$ for $\lambda$ irrational defined in the introduction is
therefore an example of a submanifold not biholomorphic to
a real-algebraic one.

We now prove the second part of Theorem \ref{fullprojwflatthm}
which we can state as follows.

\begin{thm} \label{algcodim2thm}
Let $M$ be a germ of a
real-algebraic nowhere minimal generic submanifold of
codimension 2.
Then there exists a germ of a Levi-flat real-algebraic singular hypersurface
$H$ such
that $M \subset \overline{H^*}$.
Moreover, if $M$ is not Levi-flat, then $H$ is unique.
\end{thm}

\begin{proof}%[Proof of Theorem \ref{algcodim2thm}]
If $\Orb_p$ is of constant codimension 2 in $M$, then we note that
since normal coordinates are obtained by implicit function theorem 
and there exists an algebraic implicit function theorem, then we can find
algebraic normal coordinates where $M$ is given by $w =
Q(z,\bar{z},\bar{w})$.  See \cite{BER:book} for the construction of the
normal coordinates.  Since $\Orb_p$ is of constant dimension 2 in $M$,
it agrees locally with its intrinsic complexification which is then given
by keeping $w$ constant.  Thus the vector fields
$\frac{\partial}{\partial z_k}$ and
$\frac{\partial}{\partial \bar{z}_k}$ for all $1 \leq k \leq n$ annihilate
the defining equations for $M$ (on $M$ and since $M$ is generic, in a
neighbourhood).  Thus $M$ is given by 
$w_1 = Q_1(\bar{w}_1,\bar{w_2})$ and
$w_2 = Q_2(\bar{w}_1,\bar{w_2})$.  From this we can easily construct two
algebraic holomorphic functions which are real valued on $M$, and we are done.

So assume that $\Orb_p$ is of codimension 1 in $M$ on an open and dense set.
Fix a certain representative of the germ of $M$.
Pick a point $p \in M$ near the origin
where $\Orb_p$ is of constant dimension 1 in $M$.
Let $U$ be a suitable neighbourhood of $p$.  And let $p \in U' \subset U$
be a smaller neighbourhood such that
if $\fS_{k}(q,U)$ is the $k$th Segre manifold at $q \in U'$,
$\fS_{k}(q,U)$ is connected.
We'll call $\sU$
the ambient space of $\fS_{k}(q,U)$, that is the
$U \times {}^*U \times U \times \ldots \times {}^*U\times U$
or
$U \times {}^*U \times U \times \ldots \times U \times {}^*U$
depending on whether $k$ is even or odd.  Then again
denote by $\pi \colon \C^N \times \ldots \times \C^N \to \C^N$
the projection onto the first factor, but we will define
$\pi$ on the space $\sU \times U'$.  Then define
\begin{equation*}
\fS_k(M\cap U',U) := \{ (\chi , q) \in \sU \times U' \mid
                 \chi \in \fS_k(q,U) , q \in M \cap U' \} .
\end{equation*}
$\fS_k(M\cap U',U)$ is a real-algebraic set in $\sU \times U'$ and
thus $\pi(\fS_k(M \cap U',U))$ is semialgebraic by
Tarski-Seidenberg.  We know that if $U$ is small enough and $k$ is large
enough then $\sX_q$ will lie in $\pi(\fS_k(M \cap U',U))$, and further
they form a nonsingular Levi-flat hypersurface at that point.  Since a 
semialgebraic set is contained in an algebraic set of the same
dimension, that is, there exists a polynomial $p$
defining a hypersurface $H = \{ \xi \in \C^N \mid p(\xi,\bar{\xi}) = 0 \}$
that contains $\pi(\fS_k(M \cap U',U))$.  Since $\pi(\fS_k(M \cap U',U))$
locally agrees with
a nonsingular Levi-flat hypersurface we can take $H$ to be irreducible.  Then
$H$ is Levi-flat at $p$ and by Lemma \ref{flatsomewhereflateverywhere}
it is Levi-flat.

As germs at $p$ we can see that $M \subset \overline{H^*}$.
Further, since this happens at every point where $\Orb_p$ is of
codimension 1 in $M$, and these are open and dense in $M$, then 
this must happen in some neighbourhood of the origin and hence as
germs at the origin.  Uniqueness was proved previously already
in \S \ref{sec:flatness}.
\end{proof}

\section{Example} \label{sec:ex}

%\enlargethispage{\baselineskip}

Let $M_\lambda$, $\lambda \in \R$, be the generic, nowhere minimal
submanifold of $\C^3$, with
holomorphic coordinates $(z,w_1,w_2)$ defined by
\begin{align*}
\bar{w}_1 &= e^{iz\bar{z}} w_1 ,
\\
\bar{w}_2 &= e^{i\lambda z\bar{z}} w_2 .
\end{align*}

We wish to classify the $\lambda$'s for which $M_\lambda$ has the
modulus uniqueness property at the origin.
That is, we will wish to find out when there exists a nontrivial
meromorphic function which is real valued on
$M_\lambda$.  Note that we can always find a multivalued function which is
real valued on $M_\lambda$, and that is
\begin{equation*}
(z,w_1,w_2) \mapsto
\frac{w_1^\lambda}{w_2} .
\end{equation*}
In fact, this proves that $M_\lambda$ is nowhere minimal.
Further, if $\lambda$ is rational, say $\lambda = a/b$, then
$(z,w_1,w_2) \mapsto
w_1^a / w_2^b$ is a meromorphic function that is real valued on $M$.
Thus $M_\lambda$ does not have the modulus uniqueness property, and further,
since it is of codimension 2, it is not almost minimal at the origin.

Let's check that $M_\lambda$ is almost minimal at 0 when
$\lambda$ is irrational.
For this we need to compute the Segre sets.
We can compute the third Segre set at $(z^0,w_1^0,w_2^0)$,
where $w_1^0 \not= 0$ and $w_2^0 \not= 0$,
by the following mapping (see \cite{BER:book})
\begin{equation*}
(t_1,t_2,t_3) \mapsto
(t_3,
\overline{w_1^0} e^{i(t_3 t_2 - t_2 t_1 + t_1 \overline{z^0})},
\overline{w_2^0} e^{i\lambda(t_3 t_2 - t_2 t_1 + t_1 \overline{z^0})} ) .
\end{equation*}
We can pick $t_3$ to be anything we want, and we can pick $t_2$ and $t_1$
such that the second component is anything we want since $w_1^0$ is non
zero.  By adding multiples of $2\pi$, we can add a dense set of rotations
of the third component because $\lambda$ is irrational.
This means, that the closure of this set will be
5 dimensional, and thus we will not be able to fit it inside a proper complex
analytic subset and so $M_\lambda$ is almost minimal.

We give an alternative more direct proof that $M_\lambda$ does not have
the modulus uniqueness property at the origin,
and in fact prove a slightly more general theorem
that can be used for generating further examples.

\begin{prop}
Suppose that $M$ is a real-analytic, generic submanifold of
codimension $d$ inside $\C^{n+d}$
passing through the origin
that can be defined by normal coordinates of the form
\begin{equation*}
w_j = Q_j(z,\bar{z}) \bar{w}_j
\end{equation*}
and further suppose that for any integer $K$ the functions
$Q_1^{k_1}\cdot Q_2^{k_2} \cdot \ldots \cdot Q_d^{k_d}$ for  
$0 \leq k_1 , \ldots,  k_d \leq K$ are linearly independent
as functions.  Then
there does not exist a non-constant meromorphic (nor a holomorphic) function $h$
defined in a neighbourhood of 0 which is real valued on $M$.
\end{prop}

\begin{proof}
For easier notation we will assume $n=1$ and $d=2$.
So suppose that $h = f/g$ is real valued on $M$, meaning
that on $M$ we have
$f\bar{g}-\bar{f}g = 0$.
We have proved before
that $h$ does not depend on $z$.
Suppose that $f$ and $g$ are defined by Taylor series expansions about
0.  Thus
\begin{align*}
f(w_1,w_2)
& =
\sum_{k,l \geq 0}
f_{kl} w_1^k w_2^l ,
\\
g(w_1,w_2)
& =
\sum_{n,p \geq 0}
g_{np} w_1^n w_2^p .
\end{align*}
On $M$ we therefore have (as $\bar{w}_i = \bar{Q}_i w_i)$
\begin{equation*}
\begin{split}
0 &= f\bar{g}-\bar{f}g
\\
& = 
 \left(
  \sum_{k,l,n,p \geq 0}
  f_{kl} \bar{g}_{np} w_1^k w_2^l \bar{w}_1^n \bar{w}_2^p
 \right)
 -
 \left(
  \sum_{k,l,n,p \geq 0}
  \bar{f}_{kl} g_{np} w_1^n w_2^p \bar{w}_1^k \bar{w}_2^l
 \right)
\\
& = 
\sum_{s,t \geq 0}
\left(
\sum_{k+n = s, ~ l+p = t}
  f_{kl} \bar{g}_{np} (\bar{Q}_1)^p (\bar{Q}_2)^n
 -
  \bar{f}_{kl} g_{np} (\bar{Q}_1)^l (\bar{Q}_2)^k
\right)
w_1^s w_2^t .
\end{split}
\end{equation*}

For a fixed $z$ then we have a holomorphic function in $w_1$ and $w_2$
that is 0 on a generic manifold (restriction of $M_\lambda$ to the $(w_1,w_2)$
space) and is thus identically zero.  This means that each coefficient is 0,
and since by assumption these are linear combinations of powers of
$Q_1$ and $Q_2$, we get
\begin{equation*}
  f_{(s-k)(t-l)} \bar{g}_{kl}
 -
  \bar{f}_{kl} g_{(s-k)(t-l)}
= 0 .
\end{equation*}
The above is true for all $s,t \geq 0$ and all $k \leq s, l \leq t$.
This implies that either $g \equiv 0$ or that
$f_{kl} = C g_{kl}$ for all $k,l$ for some constant $C$.  Meaning there
is no nonconstant meromorphic function which is real valued on $M$.
\end{proof}

\enlargethispage{\baselineskip}

%FIXME: else I don't get links, weird
\renewcommand{\MRhref}[2]{%
  \href{http://www.ams.org/mathscinet-getitem?mr=#1}{#2}
}

\bibliographystyle{amsalpha}
\bibliography{inftypemfld}

\providecommand{\bysame}{\leavevmode\hbox to3em{\hrulefill}\thinspace}
\providecommand{\MR}{\relax\ifhmode\unskip\space\fi MR }
\providecommand{\MRhref}[2]{%
  \href{http://www.ams.org/mathscinet-getitem?mr=#1}{#2}
}
\providecommand{\href}[2]{#2}
\begin{thebibliography}{BER00}

\bibitem[Bed77]{bedford:flat}
Eric Bedford, \emph{Holomorphic continuation of smooth functions over
  {L}evi-flat hypersurfaces}, Trans. Amer. Math. Soc. \textbf{232} (1977),
  323--341. \MR{58:1246}

\bibitem[BER98]{BER:craut}
M.~Salah Baouendi, Peter Ebenfelt, and Linda~Preiss Rothschild, \emph{C{R}
  automorphisms of real analytic manifolds in complex space}, Comm. Anal. Geom.
  \textbf{6} (1998), no.~2, 291--315. \MR{99i:32024}

\bibitem[BER99]{BER:book}
\bysame, \emph{Real submanifolds in complex space and their mappings},
  Princeton Mathematical Series, vol.~47, Princeton University Press,
  Princeton, NJ, 1999. \MR{2000b:32066}

\bibitem[BER00]{BER:localprop}
\bysame, \emph{Local geometric properties of real submanifolds in complex
  space}, Bull. Amer. Math. Soc. (N.S.) \textbf{37} (2000), no.~3, 309--336.
  \MR{2001a:32043}

\bibitem[BG99]{burnsgong:flat}
Daniel Burns and Xianghong Gong, \emph{Singular {L}evi-flat real analytic
  hypersurfaces}, Amer. J. Math. \textbf{121} (1999), no.~1, 23--53.
  \MR{2000j:32062}

\bibitem[BM88]{BM:semisub}
Edward Bierstone and Pierre~D. Milman, \emph{Semianalytic and subanalytic
  sets}, Inst. Hautes \'Etudes Sci. Publ. Math. (1988), no.~67, 5--42.
  \MR{89k:32011}

\bibitem[HJY01]{HJY:example}
Xiaojun Huang, Shanyu Ji, and Stephen S.~T. Yau, \emph{An example of a real
  analytic strongly pseudoconvex hypersurface which is not holomorphically
  equivalent to any algebraic hypersurface}, Ark. Mat. \textbf{39} (2001),
  no.~1, 75--93. \MR{2001m:32070}

\bibitem[{\L}oj91]{LS:book}
Stanis{\l}aw {\L}ojasiewicz, \emph{Introduction to complex analytic geometry},
  Birkh\"auser Verlag, Basel, 1991, Translated from Polish by Maciej Klimek.
  \MR{92g:32002}

\bibitem[Sta96]{Stanton:crauthyp}
Nancy~K. Stanton, \emph{Infinitesimal {CR} automorphisms of real
  hypersurfaces}, Amer. J. Math. \textbf{118} (1996), no.~1, 209--233.
  \MR{97h:32027}

\bibitem[Whi72]{whitney:cav}
Hassler Whitney, \emph{Complex analytic varieties}, Addison-Wesley Publishing
  Co., Reading, Mass.-London-Don Mills, Ont., 1972. \MR{52:8473}

\end{thebibliography}

\end{document}